\documentclass[12pt]{article}
\usepackage{bbm}
\usepackage{latexsym}
\usepackage{mathrsfs}
\usepackage{graphicx} %double-figure
\usepackage{subfigure}
\usepackage{amssymb}
\usepackage{amsmath}
\usepackage{amsthm}
\usepackage{color}
\usepackage{cite}
\usepackage{float}
\usepackage{indentfirst}
\usepackage{anysize}\marginsize{45mm}{45mm}{40mm}{50mm}
\usepackage{caption}
\usepackage{float}
\usepackage{multirow}

\usepackage{enumerate}

\usepackage{latexsym, bm}
\newtheoremstyle{lemma}{\topsep}{\topsep}%
     {}%         Body font
     {}%         Indent amount (empty = no indent, \parindent = para indent)
     {\bfseries}% Thm head font
     {}%        Punctuation after thm head
     {0.1em}%     Space after thm head (\newline = linebreak)
     {\thmname{#1}\thmnumber{ #2}\thmnote{ #3}}%         Thm head spec
\theoremstyle{lemma}  %%定理引用参考文献格式

\newtheorem{theorem}{Theorem}     %全局定义
\newtheorem{lemma}[theorem]{Lemma}

\newtheorem{proposition}[theorem]{Proposition}

\numberwithin{equation}{section}

\begin{document}

\title{Fractional matching preclusion for restricted hypercube-like graphs\thanks{This research was supported by the NSFC (Nos. 11801061 and 11761056), the Ministry of Education Chunhui Project (No. Z2017047) and the fundamental research funds for the central universities (No. ZYGX2018J083)}}

\author{ Huazhong L\"{u}$^{1,}$\thanks{
\emph{E-mail address}: lvhz08@lzu.edu.cn} and Tingzeng Wu$^{2}$\\
{\small $^{1}$School of Mathematical Sciences, University of Electronic Science and Technology of China,} \\
{\small  Chengdu 610054, P.R. China}\\
{\small $^{2}$School of Mathematics and Statistics, Qinghai Nationalities University, }\\
{\small Xining, Qinghai 810007, P.R. China} \\}

\date{}
\maketitle
\begin{abstract}
The restricted hypercube-like graphs, variants of the hypercube, were proposed as desired interconnection networks of parallel systems. The matching preclusion number of a graph is the minimum number of edges whose deletion results in the graph with neither perfect matchings nor almost perfect matchings. The fractional perfect matching preclusion and fractional strong perfect matching preclusion are generalizations of the concept matching preclusion. In this paper, we obtain fractional matching preclusion number and fractional strong matching preclusion numbers of restricted hypercube-like graphs, which extend some known results.
\vskip 0.1 in

\noindent \textbf{Key words:} Interconnection networks; Fractional matching preclusion; Fractional strong matching preclusion; Restricted hypercube-like graphs
\end{abstract}

\section{Introduction}
The underlying network plays important role in parallel systems.
The $n$-dimension hypercube (or binary $n$-cube), written as $Q_n$, is a well-known topology in parallel computing. To achieve desired performance that the hypercube does not have, numerous variants of the hypercube have been proposed. One among them, the hypercube-like graph, was proposed by Vaidya \cite{Vaidya} in 1993. It has been attracted considerable attention due to its outstanding performance. For example, some embedding properties, especially Hamiltonian cycle and path embeddings of the restricted hypercube-like were studied in \cite{Dong,Hsieh,Kim,Park00}. The matching preclusion number of the restricted hypercube-like graphs were determined in \cite{Park0}.

A {\em matching} is a function $f$ that each edge of $G$ is assigned a number in $\{0,1\}$ so that $\sum_{e\sim v}f(e)\leq1$ for each vertex $v\in V(G)$, where $e\sim v$ means that the sum is taken over all edges incident to $v$. A matching is {\em perfect} if $\sum_{e\sim v}f(e)=1$ for each vertex $v$, so $\sum_{e\in E(G)}=\frac{|V(G)|}{2}$. A matching is {\em almost perfect} if there exists exactly one vertex $u$ such that $\sum_{e\sim u}f(e)=0$ and $\sum_{e\sim v}f(e)=1$ for each vertex $v\in V(G)\setminus\{u\}$, so $\sum_{e\in E(G)}f(e)=\frac{|V(G)|-1}{2}$. A {\em fractional matching} is a function $f$ that each edge of $G$ is assigned a number in $[0,1]$ so that $\sum_{e\sim v}f(e)\leq1$ for each vertex $v\in V(G)$, so $\sum_{e\in E(G)}f(e)\leq\frac{|V(G)|}{2}$. Clearly, if $f(e)\in\{0,1\}$ for each edge $e$, then $f$ is a matching. If a fractional matching $f$ satisfy $\sum_{e\sim v}f(e)=1$ for each vertex $v\in V(G)$, then $f$ is a {\em fractional perfect matching} of $G$.

For $F\subseteq E(G)$, if $G-F$ has no perfect matching in $G$, then $F$ is called a {\em matching preclusion set} of $G$. The {\em matching preclusion number}, denoted by mp$(G)$, is defined to be the minimum cardinality among all matching preclusion sets. Any such set of size mp$(G)$ is called an {\em optimal matching preclusion set} (or optimal solution). This concept was proposed by Brigham et al. \cite{Brigham} as a measure of robustness of networks in the event of edge failure, as well as a theoretical connection with conditional connectivity and ``changing and unchanging of invariants''. Therefore, networks of larger mp$(G)$ signify higher fault tolerance under edge failure assumption. It is obvious that the edges incident to a common vertex form a matching preclusion set. Any such set is called a {\em trivial solution}. Therefore, mp$(G)$ is no greater than $\delta(G)$. A graph is {\em super matched} if mp$(G)=\delta(G)$ and each optimal solution is trivial. In 2011, Park et al. \cite{Park} generalized the concept of matching preclusion to strong matching preclusion as follows. A set $F$ of edges and vertices of $G$ is called a {\em strong matching preclusion set} (SMP set for short) if $G-F$ has neither perfect matching nor almost perfect matching. The {\em strong matching preclusion number} (SMP number for short) of $G$, denoted by $smp(G)$, is the minimum size of all SMP sets of $G$. The (strong) matching preclusion number of many famous interconnection networks have been investigated in the literature \cite{Cheng,Cheng1,Cheng2,Cheng3,Cheng4,Hu,Lu,Mao}

Recently, Liu and Liu \cite{Liu} generalized matching preclusion and strong matching preclusion by precluding fractional perfect matching in graphs. A set $F$ of edges of $G$ is called a {\em fractional matching preclusion set} (FMP set for short) if $G-F$ has no perfect matchings. The {\em fractional matching preclusion number} (FMP number for short) of $G$, denoted by $fmp(G)$, is the minimum size of all FMP sets of $G$. Clearly, $fmp(G)\leq \delta(G)$. Moreover, by the definition of $fmp(G)$, if $G$ has even order, then $\delta(G)\leq fmp(G)$. A set $F$ of edges and vertices of $G$ is called a {\em fractional strong matching preclusion set} (FSMP set for short) if $G-F$ has no fractional perfect matchings. The {\em fractional strong matching preclusion number} (FSMP number for short) of $G$, denoted by $fsmp(G)$, is the minimum size of all FSMP sets of $G$.

The fractional perfect (strong) matching preclusion number of $(n,k)$-star graphs has been determined in \cite{Ma}. In \cite{Liu}, the authors obtained fractional perfect (strong) matching preclusion number the complete graph, the Petersen graph and the twisted cube. In this paper, we determine fractional perfect (strong) matching preclusion number of restricted hypercube-like graphs, which include the twisted cubes as a proper subset.

The rest of this paper is organized as follows. In Section 2, some notations, the definitions of the balanced hypercube and some useful lemmas are presented. Section 3 shows the existence of two edge-disjoint Hamiltonian cycles of the balanced hypercube and provides an algorithm to construct two edge-disjoint Hamiltonian cycles of the balanced hypercube. Finally, conclusions are given in Section 4.

\section{Preliminaries}

In this section, we shall present some notations, definitions of the restricted hypercube-like graphs and some useful lemmas.

\vskip 0.1 in

Interconnection networks are usually modeled by graphs, where vertices represent processors and edges represent links between processors. Throughout this paper, we only consider finite and simple undirected graphs. Let $G=(V(G),E(G))$ be a graph, where $V(G)$ is the vertex-set of $G$ and $E(G)$ is the edge-set of $G$. The number of vertices of $G$ is denoted by $|V(G)|$. Two
vertices $u$ and $v$ are adjacent if $uv\in E(G)$. A {\em neighbor} of a vertex $v$ in $G$ is any vertex incident to $v$. A {\em path} $P$ in $G$ is a sequence of distinct vertices so that there is an edge joining each pair of consecutive vertices. If $P=v_0v_1\cdots v_{k-1}$ is a path and $k\geq3$, then the graph $C=P+v_{k-1}v_0$ is said to be a {\em cycle}. The above path $P$ and cycle $C$ might be written as $\langle v_0,v_1,\cdots ,v_{k-1}\rangle$ and $\langle v_0,v_1,\cdots ,v_{k-1},v_0\rangle$, respectively. The length of a path or a cycle is its number of edges. A cycle of length $k$ is called a $k$-{\em cycle}. A cycle containing all vertices of a graph $G$ is called a {\em Hamiltonian cycle}. A graph is called {\em $f$-fault Hamiltonian} (resp. {\em $f$-fault Hamiltonian-connected)} if there exists a Hamiltonian cycle (resp. if there exists a Hamiltonian path joining each pair of vertices) in $G-F$ for any set $F$ of vertices and/or edges with $|F|\leq f$. For other standard graph notations and terminologies not defined here please refer to \cite{Bondy}.

\vskip 0.1 in

Let $G_0$ and $G_1$ be two disjoint graphs with the same order. In addition, let $\Phi(G_0,G_1)$ be all bijections from $V(G_0)$ to $V(G_1)$. Given a bijection $\phi\in\Phi(G_0,G_1)$, $v\phi(v)$ is an edge from $G_0$ to $G_1$ for any $v\in V(G_0)$, and we denote $G_0\oplus_{\phi} G_1$ a graph whose vertex set is $V(G_0)\cup V(G_1)$ and edge set is $E(G_0)\cup E(G_1)\cup \{v\phi(v)|v\in V(G_0)\}$. Clearly, each vertex in $G_0$ has exact one neighbor in $G_1$, and vice versa. When the context is clear, we often omit the symbol $\phi$ from $\oplus_{\phi}$. By using the above graph operator, Vaidya et al. \cite{Vaidya} gave a recursive definition of the hypercube-like graphs as follows. $HL_0=\{K_1\}$ and $HL_m=\{G_0\oplus_{\phi}G_1|G_0,G_1\in HL_{m-1},\phi\in\Phi(G_0,G_1)\}$ for $m\geq1$. A graph $HL_m$ is called {\em $m$-dimensional HL-graph}. The {\em restricted HL-graphs}, which is an interesting subset of HL-graphs, were proposed by Park et al. in \cite{Park2}. $RHL_3=HL_3\setminus Q_3=\{G(8,4)\}$, $RHL_m=\{G_0\oplus_{\phi}G_1|G_0,G_1\in RHL_{m-1},\phi\in\Phi(G_0,G_1)\}$ for $m\geq4$, where $Q_3$ is the 3-dimensional hypercube, and $G(8,4)$ is the recursive circulant whose vertex set is $\{v_i|0\leq i\leq 7\}$ and edge set is $\{v_iv_j|j\equiv i+1$ or $i+4($mod $8)\}$. Any graph contained in $RHL_m$ is called an {\em $m$-dimensional restricted HL-graph} and is denoted by $G^m$. Since $G(8,4)$ is nonbipartite, the restricted HL-graphs are all nonbipartite, forming a proper subset of nonbipartite HL-graphs. It is noticeable that numerous of famous interconnection networks such as crossed cube, M\"{o}bius cube, twisted cube, Mcube, generalized twisted cube are known to be restricted HL-graphs \cite{Park2}.

\vskip 0.1 in

In what follows, we shall present some useful results.

\begin{proposition}\label{fpm}\cite{Scheinerman}{\bf.} A graph $G$ has a
fractional perfect matching if and only if $i(G-S)\leq|S|$ for every set $S\subseteq V(G)$,
where $i(G-S)$ is the number of isolated vertices of $G-S$.
\end{proposition}

\begin{lemma}\label{fault-Hamiltonian}\cite{Park2}{\bf .} Every $G^m$ with $m\geq3$ is ($m-3$)-fault Hamiltonian-connected and ($m-2$)-fault Hamiltonian.
\end{lemma}

\begin{theorem}\label{smp-Gm}\cite{Park}{\bf .} For each $m\geq3$, $smp(G^m)=m$.
\end{theorem}

\section{Main results}

Since all restricted hypercube-like graphs have even order, combining Theorem \ref{smp-Gm}, the following theorem is obvious.

\begin{theorem}{\bf.}
For each integer $m\geq3$, $fmp(G^m)=m$.
\end{theorem} 

Next we consider FSMP number of $G^m$. The following lemma gives both lower and upper bounds of FSMP number of $G^m$.

\begin{lemma}\label{Gm-fsmp-bound}{\bf.} Let $m\geq3$ be an integer. Then $m-1\leq fsmp(G^m)\leq m$.
\end{lemma}
\noindent{\bf Proof.} Since $G^m$ is $m$-regular, $fsmp(G^m)\leq m$. It suffices to show that there exist no FSMP sets of $G^m$ with size at most $m-2$. Let $F\subseteq V(G^m)\cup E(G^m)$ with $|F|\leq m-2$. By Lemma \ref{fault-Hamiltonian}, $G^m-F$ has a Hamiltonian cycle $C$. Assigning $\frac{1}{2}$ to each edge of $C$ and 0 to other edges, we can obtain a fractional perfect matching of $G^m-F$. Hence, $fsmp(G^m)\geq m-1$. This completes the proof.\qed

To obtain the exact value of $fsmp(G^m)$, we begin with $m=3$.

\begin{lemma}\label{G3-fsmp}{\bf.} $fsmp(G^3)=2$. Additionally, the optimal FSMP set contains exactly one vertex $u$ and one boundary edge $e$, where $u$ is adjacent to one of the end vertices of $e$.
\end{lemma}

\begin{figure}
\centering
\includegraphics[width=100mm]{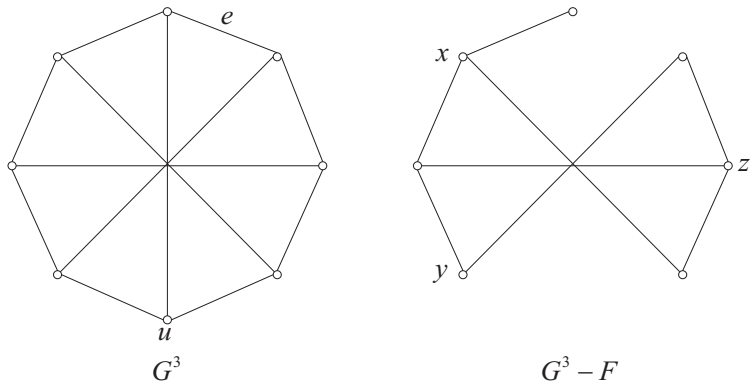}
\caption{An FSMP set of $G^3$.} \label{g2}
\end{figure}

\noindent{\bf Proof.} By Lemma \ref{Gm-fsmp-bound}, we have $fsmp(G^3)\geq2$. It can be verified that if $e$ is a diagonal edge of $G^3$, then $(G^3-e)-v$ contains a fractional perfect matching for any vertex $v\in V(G^3-e)$. Thus, by symmetry of $G^3$, let $e$ be any boundary edge of $G^3$. Observe that there are exact two vertices of $G^3$ such that they are nonadjacent to any end vertices of $e$. By deleting exact one of them from $G^3-e$, it is easy to find a fractional perfect matching of the resulting graph. Observe also that there are exact four vertices of $G^3$ such that they are adjacent to one of end vertices of $e$. Without loss of generality, suppose that $u$ is such a vertex of $G^3$. Then we can find three vertices $x,y$ and $z$ of $G^3$ such that $i(((G^3-e)-u)-\{x,y,z\})=4>3$ (see Fig. \ref{g2}). By Proposition \ref{fpm}, $(G^3-e)-u$ has no fractional perfect matchings. Thus, $F$ is an FSMP set of $G^3$. Accordingly, $fsmp(G^3)=2$. Moreover, since $mp(G^3)=smp(G^3)=3$, any optimal FSMP set must contain one vertex and one edge. This completes the proof.\qed

\begin{figure}
\centering
\includegraphics[width=50mm]{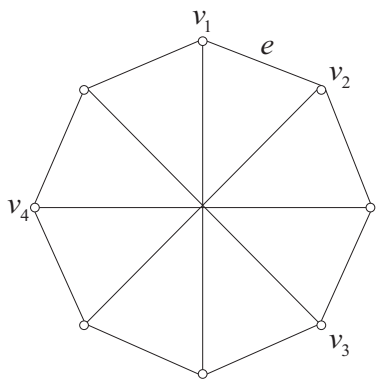}
\caption{A remainder set of $G^3$.} \label{g3}
\end{figure}

As stated in the proof of Lemma \ref{G3-fsmp}, if an optimal FSMP set $F$ is deleted from $G^3$, the resulting graph contains three vertices $x, y$ and $z$ such that $i(G^3-F-\{x,y,z\})=4$. So we denote the four isolated vertices in $G^3-F-\{x,y,z\}$ by $v_1, v_2, v_3$ and $v_4$, respectively. Observe that we have two essentially different choices of $F$ by symmetry, nevertheless, the position of $v_1, v_2, v_3$, and $v_4$ on $G^3$ is unique under isomorphism (see Fig. \ref{g3}). For convenience, any set of four vertices with the same position as $v_1, v_2, v_3$, and $v_4$ is called a {\em remainder set} of $G^3$. Since $G^4$ is constructed from two copies of $G^3$ by joining a perfect matching. In the following lemma, we characterize $fsmp(G^4)$.

\begin{lemma}\label{G4-fsmp-necessity}{\bf.} $fsmp(G^4)=3$ If there exists a remainder set $R$ of $G_0$ such that the resulting graph of $G_1-\phi(R)$ contains at most one edge, then $fsmp(G^4)=3$. Otherwise, $fsmp(G^4)=4$.
\end{lemma}

\noindent{\bf Proof.} By Lemma \ref{Gm-fsmp-bound}, $3\leq fsmp(G^4)\leq 4$. Let $F\subseteq V(G^4)\cup E(G^4)$ with $|F|=3$ and let $F_V=F\cap V(G^4)$ and $F_E=F\cap E(G^4)$. We may assume that $|F_0|=\max\{|F_0|,|F_1|\}$. We shall prove that: (1) If there exists a vertex set $S$ such that $i(G^4-F-S)>|S|$, then $fsmp(G^4)=3$; (2) If $G^4-F$ has a fractional perfect matching, then $fsmp(G^4)=4$. If $|F_V|$ is even, then $G^4-F$ has a perfect matching by Lemma \ref{smp-Gm}. It remains to consider the case that $|F_V|=1$ or $|F_V|=3$. We distinguish the following cases.

\noindent{\bf Case 1.} $|F_0|=1$. Obviously, $|F_1|\leq1$. By Lemma \ref{G3-fsmp}, each of $G_i$, $i=0,1$, has a fractional perfect matching $f_i$. Then a fractional perfect matching $f$ of $G^4-F$ can be obtained as follows.

$$ f(e)=\left\{
\begin{aligned}
f_i(e), & \ \ \ e\in E(G_i-F_i); \\
0, &  \ \ \ e\in E_c-F_c.
\end{aligned}
\right.
$$

\noindent{\bf Case 2.} $|F_0|=2$. Since $|F|=3$, $|F_1|\leq1$. By Lemma \ref{G3-fsmp}, if $F_0$ contains two vertices or two edges, respectively, then $G_0-F_0$ has a perfect matching. Similar to the proof of Case 1, $G^4-F$ has a fractional perfect matching. So we assume that $F_0$ contains exact one vertex and exact one edge. Let $F_0=\{u,e\}$, where $u\in V(G^4)$ and $e\in E(G^4)$. If $G_0-F_0$ has a fractional perfect matching, again, $G^4-F$ has a fractional perfect matching, we are done. So we assume that $G_0-F_0$ contains no fractional perfect matchings. We further distinguish the following subcases.

%

%\noindent{\bf Subcase 2.1.} $F_1\subset F_V$. Let $F_1=\{x\}$. Since each vertex of $G_0$ has exact one neighbor in $G_1$, we may assume that $v$ has a neighbor $v'$ ($v'\neq x$) in $G_1-F_1$. By Lemma \ref{G3-fsmp}, $G_1-F_1-v'$ has a fractional perfect matching $f_1$. Let $f_0$ be a fractional perfect matching of $G_0-F_0-v$. Then a fractional perfect matching $f$ of $G^4-F$ can be obtained as follows.
%
%$$ f(e)=\left\{
%\begin{aligned}
%f_0(e), & \ \ \ e\in E(G_0-F_0-v); \\
%f_1(e), & \ \ \ e\in E(G_1-F_1-v'); \\
%1, &  \ \ \ e=vv'; \\
%0, &  \ \ \ e\in E_C-F_C.
%\end{aligned}
%\right.
%$$

\noindent{\bf Subcase 2.1.} $F_1\subset F_E$. By Proposition \ref{fpm}, there exists a set of vertices $S_0=\{x,y,z\}$ such that $i(G_0-F_0-S_0)=4$. Correspondingly, $G_0-F_0-S_0$ contains exact four isolated vertices. Without cause of ambiguity, let $R=G_0-F_0-S_0$ be a remainder set of $G_0$. Let $\phi(R)\subset V(G_1)$ be the set containing all neighbors of vertices of $R$ in $G_1$. For any edge $e'\in E(G_1)$, $G_1-e'$ is Hamiltonian by Lemma \ref{fault-Hamiltonian}. Thus, for any $S_1\subset V(G_1)$, $i(G_1-e'-S_1)\leq |S_1|$ holds. If $G^4-F$ has no fractional perfect matchings, then $|S_1|=4$. In fact, $S_1=\phi(R)$. If not. Then there exists an edge joining a vertex in $R$ and a vertex in $G_1-S_1$, which implies that $i(G^4-F-S_0\cup S_1)<7$, a contradiction. Moreover, we claim that $G_1-S_1$ contains at most one edge. On the contrary, assume that $G_1-S_1$ contains at least two edges. We have $i(G^4-F-S_0\cup S_1)<7$, a contradiction again. It implies that if $R$ is a remainder set of $G_0$ and $G_1-\phi(R)$ contains at most one edge, then $fsmp(G^4)=3$; Otherwise, $fsmp(G^4)=4$.

\noindent{\bf Subcase 2.2.} $F_1=\emptyset$. Obviously, $|F_C|=1$. By Lemma \ref{fault-Hamiltonian}, $G_0-F_0$ has either a Hamiltonian cycle or a Hamiltonian path with seven vertices. If $G_0-F_0$ has a Hamiltonian cycle, then there is a fractional perfect matching of $G_0-F_0$. Observe that $G_1-F_1$ has a fractional perfect matching, we have that $G^4-F$ has a fractional perfect matching. So we assume that $G_0-F_0$ has a Hamiltonian path $P$ with seven vertices. Let $v$ and $w$ be end vertices of $P$. Since each vertex of $G_0$ is incident to a cross edge, $v$ or $w$, say $v$, has a neighbor $v'\in V(G_1)$ such that $vv'\not\in F_C$. By Lemma \ref{G3-fsmp}, $G_1-v'$ has a fractional perfect matching $f_1$. Let $f_0$ be a fractional perfect matching of $G_0-F_0-v$. Then a fractional perfect matching $f$ of $G^4-F$ can be obtained as follows.

$$ f(e)=\left\{
\begin{aligned}
f_0(e), & \ \ \ e\in E(G_0-F_0-v); \\
f_1(e), & \ \ \ e\in E(G_1-v'); \\
1, &  \ \ \ e=vv'; \\
0, &  \ \ \ e\in E_C-F_C.
\end{aligned}
\right.
$$

\noindent{\bf Case 3.} $|F_0|=3$. By our assumption, $F_0$ contains exact one vertex or three vertices. By Lemma \ref{fault-Hamiltonian}, $G_0-F_0$ has a Hamiltonian cycle, or has a Hamiltonian path with odd number of vertices, or has a spanning subgraph containing one odd path and one even path. Similar to the proof above, we only consider the case that $G_0-F_0$ contains one odd path $P$ and one even path $Q$. We assume that $u$ is one of the end vertices of the odd path and $uu'$ is the cross edge. Obviously, $Q$ has a perfect matching $f_Q$ and $P-u$ has a perfect matching $f_{P-u}$. Moreover, $G_1-u'$ has a Hamiltonian cycle. Thus, $G_1-u'$ has a fractional perfect matching $f_{G_1-u'}$. Then a fractional perfect matching $f$ of $G^4-F$ can be obtained as follows.

$$ f(e)=\left\{
\begin{aligned}
f_Q(e), & \ \ \ e\in E(Q); \\
f_{P-u}(e), & \ \ \ e\in E(P-u); \\
f_{G_1-u'}(e), & \ \ \ e\in E(G_1-u'); \\
1, &  \ \ \ e=uu'; \\
0, &  \ \ \ e\in E_C-uu'.
\end{aligned}
\right.
$$

\qed

Interestingly, in the above lemma, the necessary condition for $fsmp(G^4)=3$ is also sufficient.

\begin{lemma}\label{G4-fsmp-equivalence}{\bf.} If $fsmp(G^4)=3$, then there exists a remainder set $R$ of $G_0$ such that the resulting graph of $G_1-\phi(R)$ contains at most one edge.
\end{lemma}
\noindent{\bf Proof.} Obviously, $|F_V|$ is odd. Suppose that $|F_V|=3$, by Lemma \ref{smp-Gm}, then $G^4-F$ has a perfect matching. Thus, $G^4-F$ has a fractional perfect matching. So $|F_V|=1$ and $|F_E|=2$. Without loss of generality, suppose that $F_V=\{u\}$ and $u\in V(G_0)$. We distinguish the following cases.

\noindent{\bf Case 1.} $|F_E\cap E(G_0)|=2$. By the proof of Case 3 in Lemma \ref{G4-fsmp-necessity}, we have $G^4-F$ has a fractional perfect matching, a contradiction.

\noindent{\bf Case 2.} $|F_E\cap E(G_1)|=2$. By Lemma \ref{fault-Hamiltonian}, $G_0-F_0$ has a Hamiltonian cycle and $G_1-F_1$ has a Hamiltonian path with eight vertices. So $G^4-F$ has a fractional perfect matching, a contradiction.

\noindent{\bf Case 3.} $|F_E\cap E(G_0)|\leq1$ and $|F_E\cap E(G_1)|\leq1$. If $|F_E\cap E(G_0)|=0$, then $G_0-F_0$ has a Hamiltonian cycle. Similarly, $G_1-F_1$ has a Hamiltonian cycle. Obviously, $G^4-F$ has a fractional perfect matching, a contradiction. Thus, $|F_E\cap E(G_0)|=1$. If $G_0-F_0$ has a fractional perfect matching, then $G^4-F$ has a fractional perfect matching, a contradiction. So we assume that $G_0-F_0$ has no fractional perfect matchings. Thus, there exists a remainder set $R$ of $G_0$. Moreover, if $|F_E\cap E(G_1)|=0$, then $|F_E\cap F_C|=1$. By the proof of Subcase 2.2 of Lemma \ref{G4-fsmp-necessity}, we know that $G^4-F$ has a fractional perfect matching, a contradiction. Thus, $|F_E\cap E(G_1)|=1$. Similar to the proof of Subcase 2.1 of Lemma \ref{G4-fsmp-necessity}, the statement follows.\qed

In the following, we study FSMP number of $G^5$ as our induction basis.

\begin{lemma}\label{G5-fsmp}{\bf.} $fsmp(G^5)=5$.
\end{lemma}

\noindent{\bf Proof.} By Lemma \ref{Gm-fsmp-bound}, $4\leq fsmp(G^5)\leq 5$. Let $F\subseteq V(G^5)\cup E(G^5)$ with $|F|=4$ and let $F_V=F\cap V(G^5)$ and $F_E=F\cap E(G^5)$. We may assume that $|F_0|=\max\{|F_0|,|F_1|\}$. We shall show that $G^5-F$ has a fractional perfect matching. If $|F_V|$ is even, then $G^5-F$ has a perfect matching by Lemma \ref{smp-Gm}. It remains to consider the case that $|F_V|=1$ or $|F_V|=3$. Since $G_i\cong G^4$, $i=0,1$, each of $G_i-F_i$ has a Hamiltonian cycle if $|F_0|\leq 2$. Thus, $G^5-F$ has a fractional perfect matching. So we only consider the case that $|F_0|\geq 3$. We distinguish the following two cases.

\noindent{\bf Case 1.} $|F_0|=3$. Then $|F_1|\leq1$. If $fsmp(G_0)=4$, then $G_0-F_0$ has a fractional perfect matching and thus, $G^5-F$ has a fractional perfect matching. So we assume that $fsmp(G_0)=3$ and hence, $G_0-F_0$ has no fractional perfect matchings. By Lemma \ref{G4-fsmp-equivalence}, we have $|F_V|=1$. So $G_0-F_0$ has a Hamiltonian path $P$ with odd vertices. Let $v$ and $w$ be end vertices of $P$. Since each vertex of $G_0$ is incident to a cross edge, we may assume that $vv'\not\in F_C$ is a cross edge. By Lemma \ref{G3-fsmp}, $G_1-F_1-v'$ has a fractional perfect matching $f_1$. Let $f_0$ be a fractional perfect matching of $G_0-F_0-v$. Then a fractional perfect matching $f$ of $G^5-F$ can be obtained as follows.

$$ f(e)=\left\{
\begin{aligned}
f_0(e), & \ \ \ e\in E(G_0-F_0-v); \\
f_1(e), & \ \ \ e\in E(G_1-F_1-v'); \\
1, &  \ \ \ e=vv'; \\
0, &  \ \ \ e\in E_C-vv'.
\end{aligned}
\right.
$$

\noindent{\bf Case 2.} $|F_0|=4$. Then $|F_1|=0$. By Lemma \ref{fault-Hamiltonian}, $G_0-F_0$ has a Hamiltonian cycle, or has a Hamiltonian path with odd number of vertices, or has a spanning subgraph containing one odd path and one even path. Similar to the proof of Case 3 of Lemma \ref{G4-fsmp-necessity}, we can obtain that $G^5-F$ has a fractional perfect matching. \qed

Now we are ready to present the main result of this paper.

\begin{theorem}\label{Gm-fsmp}{\bf.} $fsmp(G^m)=m$ for $m\geq5$.
\end{theorem}

\noindent{\bf Proof.} It suffices to prove that $fsmp(G^m)\neq m-1$ by Lemma \ref{Gm-fsmp-bound}. We shall prove that for any set $F$ of vertices and edges with $|F|=m-1$, $G^m-F$ has a fractional perfect matching. We proceed by induction on $m$. By Lemma \ref{G5-fsmp}, the statement holds for $m=5$. We assume that the statement holds for all integers not greater that $m-1$ with $m\geq6$. Next we consider $G^m$. We consider two cases.

\noindent{\bf Case 1.} $|F_0|\leq m-2$. By the induction hypothesis, we have $fsmp(G^{m-1})=m-1$. So each of $G_i-F_i$, $i=0,1$, has a fractional perfect matching $f_i$. Thus, $f_0\cup f_1$ is a fractional perfect matching of $G^m-F$.

\noindent{\bf Case 2.} $|F_0|=m-1$. Obviously, $F-F_0=\emptyset$. We can choose $\alpha\in F_0$ such that $F_0-\{\alpha\}$ contains even number of vertices. For convenience, let $F_0'=F_0-\{\alpha\}$. Then $|F_0'|=m-2$. So $G_1-F_0'$ has a perfect matching $f_0$ by Lemma \ref{smp-Gm}. We further consider two subcases.

\noindent{\bf Case 2.1.} $\alpha$ is a vertex. We may assume that $\alpha=v$ and $uv$ is the edge that $f_0(uv)=1$. Then the restriction of $f_0$ on $E(G_0-F_0'-\{u,v\})$ is a perfect matching of $G_0-F_0'-\{u,v\}$. Let $uu'$ be a cross edge of $G^m$. Moreover, there exists a fractional perfect matching $f_1$ of $G_1-u'$. By assigning $uu'$ with $1$ and other cross edges 0, we can obtain a fractional perfect matching of $G^m-F$.

\noindent{\bf Case 2.2.} $\alpha$ is an edge. If $f_0(\alpha)=0$, then $f_0$ is a perfect matching of $G_0-F_0$. So we assume that $f_0(\alpha)=1$. Then the restriction of $f_0$ on $E(G_0-F_0'-\{u,v\})$ is a perfect matching of $G_0-F_0'-\{u,v\}$. We may assume that $\alpha=uv$. Thus, there exists two cross edges $uu'$ and $vv'$. Clearly, there exists a fractional perfect matching $f_1$ of $G_1-\{u',v'\}$. By assigning $uu'$ and $vv'$ with $1$ and other cross edges 0, we can obtain a fractional perfect matching of $G^m-F$.\qed

\section{Conclusions}

In this paper, we obtain FMP and FSMP number of restricted hypercube-like graphs. Matching preclusion problem has been attracted much attention in the literature. Since FMP and FSMP problems are interesting generalizations of matching preclusion problem, it is meaningful to consider FMP and FSMP number of famous interconnection networks, as well as theory of FMP in general graphs.

%\vskip 0.3 in
%
%\noindent{\bf{\normalsize Acknowledgements}}
%
%\vskip 0.1 in

%The authors are grateful to the anonymous referees for their comments and constructive suggestions that greatly improved the original manuscript.


\begin{thebibliography}{99}

\bibitem{Bondy} J.A. Bondy, U.S.R. Murty, Graph theory, Springer, New York, 2007.

\bibitem{Brigham} R. Brigham, F. Harary, E. Violin, J. Yellen, Perfect-matching preclusion, Congr. Numer. 174 (2005) 185--192.

\bibitem{Cheng} E. Cheng, P. Hu, R. Jia, L. Lipt\'{a}k, B. Scholten, J, Voss, Matching preclusion and conditional matching preclusion for pancake and burnt pancake graphs, Inter. J. Parallel, Emerg. Distr. Syst. 29 (2014) 499--512.

\bibitem{Cheng1} E. Cheng, L. Lesniak, M. Lipman, L. Lipt\'{a}k, Conditional matching preclusion sets, Inform. Sci. 179 (2009) 1092--1101.

\bibitem{Cheng2} E. Cheng, L. Lipman, L. Lipt\'{a}k, Matching preclusion and conditional matching preclusion for regular interconnection networks, Discrete Appl. Math. 160 (2012) 1936--1954.

\bibitem{Cheng3} E. Cheng, S. Shah, V. Shah, D.E. Steffy, Strong matching preclusion for augmented cubes, Theor. Comput. Sci. 491 (2013) 71--77.

\bibitem{Cheng4} E. Cheng, O. Siddiqui, Strong matching preclusion of arrangement graphs, J. Interconn. Networks 16 (2016) 1650004.

\bibitem{Dong} Q. Dong, J. Zhou, Y. Fu, H. Gao, Hamiltonian connectivity of restricted hypercube-like networks under the conditional fault model, Theor. Comput. Sci. 472 (2013) 46--59.

\bibitem{Hsieh} S.-Y. Hsieh, C.-W. Lee, Pancyclicity of restricted hypercube-like networks under the conditional fault model, SIAM J. Discrete Math. 23 (4) (2010) 2100--2119.

\bibitem{Hu} X. Hu, Y. Tian, X. Liang, J. Meng, Strong matching preclusion for n-dimensional torus networks, Theor. Comput. Sci. 635 (2016) 64--73.

\bibitem{Kim} S.-Y. Kim, J.-H. Park, Many-to-many two-disjoint path covers in restricted hypercube-like graphs, Theor. Comput. Sci. 531 (2014) 26--36.

\bibitem{Liu} Y. Liu, W. Liu, Fractional matching preclusion of graphs, J. Comb. Optim. 34 (2017) 522--533.

\bibitem{Lu} H. L\"{u}, X. Li, H. Zhang, Matching preclusion for
balanced hypercubes, Theor. Comput. Sci. 465 (2012) 10--20.

\bibitem{Ma} T. Ma, Y. Mao, E. Cheng, J. Wang, Fractional Matching Preclusion for $(n,k)$-Star Graphs, Parall. Process. Lett. 28 (2018) 1850017.

\bibitem{Mao} Y. Mao, Z. Wang, E. Cheng, C. Melekian, Strong matching preclusion number of graphs, Theor. Comput. Sci. 713 (2018) 11--20.


\bibitem{Park0} J. Park, Matching preclusion problem in restricted HL-graphs and recursive circulant $G(2^m,4)$, J. KIISE 35 (2) (2008) 60--65.

\bibitem{Park00} J.-H. Park, Paired many-to-many disjoint path covers in restricted hypercube-like graphs, Theor. Comput. Sci. 634 (2016) 24--34.


\bibitem{Park} J. Park, I. Ihm, Strong matching preclusion, Theor. Comput. Sci. 412 (2011) 6409--6419.

\bibitem{Park2} J.-H. Park, H.-C. Kim, H.-S. Lim, Fault-hamiltonicity of hypercube-like interconnection networks, in: Proc. IEEE International Parallel and Distributed Processing Symposium IPDPS 2005, Denver, April 2005.

\bibitem{Scheinerman} E.R. Scheinerman, D.H. Ullman, Fractional graph theory, Wiley, New York, 1997.

\bibitem{Vaidya} A.S. Vaidya, P.S.N. Rao, S.R. Shankar, A class of hypercube-like networks, in: Proc. of the 5th Symp. on Parallel and Distributed Processing, IEEE Comput. Soc., Los Alamitos, CA, December 1993, pp. 800--803.

\end{thebibliography}
\end{document}